\documentclass[10pt, letterpaper]{amsart}
\usepackage{amsmath,amssymb,amscd,amsthm,eucal,curves,epic, graphicx}
\usepackage[all]{xy}
\usepackage[latin1]{inputenc}
\usepackage{color}
\usepackage{floatflt}

\swapnumbers

\theoremstyle{plain}
\newtheorem{theorem}{Theorem}[section]
\newtheorem{lemma}[theorem]{Lemma}
\newtheorem{proposition}[theorem]{Proposition}

\theoremstyle{definition}

\newtheorem{remark}[theorem]{Remark}

\newtheorem{nn}[theorem]{}

\newcommand{\ra}{\rightarrow}

\newcommand{\hra}{\hookrightarrow}

\newcommand{\retrmanlabel}[4]{\ensuremath \xymatrix@1{#1\ar@<0.4ex>[r]^-{#3}&#2 \ar@<0.4ex>[l]^-{#4} }}
\newcommand{\retrlabel}[2]{\ensuremath \xymatrix@1{#1\ar@<0.4ex>[r]^-{r}&#2 \ar@<0.4ex>[l]^-{s} }}
\newcommand{\retr}[2]{\ensuremath \xymatrix@1{#1\ar@<0.4ex>[r]&#2 \ar@<0.4ex>[l] }}

\renewcommand{\SS}[1][\bullet]{\mathcal{S}_{#1}}
\newcommand{\PP}[1][\bullet]{\mathcal{P}_{#1}}
\newcommand{\TT}[1][\bullet]{\mathcal{T}_{#1}}
\newcommand{\Rfd}{\mathcal{R}^{fd}}

\newcommand{\simp}[1][B]{{\mathcal S}(#1)}

\newcommand{\Ch}{{\mathcal C}h}

\newcommand{\colim}[1][]{\ifmmode\ifinner {\operatorname{colim}_{#1}}\,\else \underset{#1}{\operatorname{colim}}\, \fi\fi}
\newcommand{\hocolim}[1][]{\ifmmode\ifinner {\operatorname{hocolim}_{#1}}\,\else \underset{#1}{\operatorname{hocolim}}\,\fi\fi}
\newcommand{\holim}[1][]{\ifmmode\ifinner {\operatorname{holim}_{#1}}\,\else \underset{#1}{\operatorname{holim}}\,\fi\fi}

\newcommand{\wh}{{\rm Wh} } 

\newcommand{\ehom}{\ensuremath\chi^h}

\newcommand{\tsm}{\ensuremath\tau^{s}_{\rho}}

\newcommand{\hofib}{\mathrm{hofib} }

\newcommand{\id}{{\ensuremath{\rm id}}}

\newcommand{\RR}{\mathcal{R}}
\newcommand{\fd}{\mathbb{F}}
\newcommand{\real}{\mathbb{R}}

\title[Smooth parametrized torsion -- a manifold approach]{Smooth parametrized torsion \\  a manifold approach}
\author[B. Badzioch]{Bernard Badzioch}
\author[W. Dorabia{\l}a]{Wojciech Dorabia{\l}a}
\author[B. Williams]{Bruce Williams}

\address[]{Department of Mathematics, University at Buffalo, SUNY, Buffalo, NY}
\address[]{Department of Mathematics, Penn State Altoona, Altoona, PA }
\address[]{Department of Mathematics, University of Notre Dame, IN}

\begin{document}

\begin{abstract}
We give a construction of a torsion invariant of bundles of smooth manifolds 
which is based on the work of Dwyer, Weiss and Williams on smooth structures on
fibrations. 
\end{abstract}

\date{11/11/2007}
\maketitle


\section{Introduction}
\label{INTRO}

\begin{nn}
Recently there has been considerable interest and activity related to the construction and 
the computations of parametrized torsion invariants. The goal here is the development of 
a generalization of Reidemeister torsion -- which is a classical secondary invariant 
of CW-complexes -- to bundles of manifolds. One approach to this problem was proposed 
by Bismut and Lott \cite{Bis-Lott}. Another, using parametrized Morse functions resulted 
from the work of Igusa \cite{Igusa} and Klein \cite{Klein}.

In \cite{DWW} Dwyer, Weiss, and Williams presented yet another definition of torsion of bundles 
whose main feature is that it is described entirely in terms of algebraic topology .  As a result
their construction is quite intuitive. Given a smooth bundle $p\colon E\to B$
we can consider the Becker-Gottlieb transfer map $p^{!}\colon B\to Q(E_{+})$. If 
 $\rho\colon M\to E$ is a locally constant sheaf of $R$-modules we can construct a map 
$c_{\rho}\colon B\to K(R)$ which assigns to a point $b\in B$ the point of $K(R)$ represented 
by the singular chain complex $C_{\ast}(F_{b}, \rho|_{F_{b}})$
of the fiber of $p$ over $b$ with coefficient in $\rho$. One of the results of \cite{DWW} implies 
that there exists a map $\lambda\colon Q(E_{+})\to K(R)$ such that the diagram 
$$\label{MAIN_DIAG}
 \xymatrix{
  &Q(E_{+})\ar[d]^{\lambda} \\
B\ar[ur]^{p^{!}}\ar[r]_{c_{\rho}} & K(R)  \\
}
$$  
commutes up to a preferred homotopy. If the sheaf $\rho$ is such that all chain complexes 
$C_{\ast}(F_{b}, \rho|_{F_{b}})$ are acyclic, then the map $c_{\rho}$ is canonically homotopic
to the constant map. Thus we obtain a lift of $p^{!}$ to the space
$\wh_{\rho}(E)$ which is  the homotopy fiber of $\lambda$ over the trivial element of $K(R)$. 
This lift $\tau_{\rho}(p)$ is the smooth torsion of the bundle $p$. 

In order to see what kind of information about a bundle is carried by its torsion lets assume 
first that we are given two smooth bundles $p_{i}\colon E_{i}\to B$, $i=1, 2$ and a map of bundles
$f\colon E_{1}\to E_{2}$. Let $A(E_{i})$ denote the Waldhausen $A$-theory of the space $E_{i}$. 
We have the assembly maps $a\colon Q(E_{i+})\to A(E_{i})$ (see Section \ref{ASSEMBLY}) which
fit into a commutative square
$$
\xymatrix{
Q(E_{1+}) \ar[r]^{f_{\ast}}\ar[d]_{a} & Q(E_{2+})\ar[d]^{a} \cr
A(E_{1}) \ar[r]_{f_{\ast}}  & A(E_{2}) \cr
}
$$
Let $p^{A}_{i}\colon B\to A(E_{i})$ be the composition $p^{A}_{i}=a p_{i}^{!}$.
One of the implications of \cite{DWW} is that if $f$ is a fiberwise homotopy equivalence 
then we can construct a homotopy $\omega_{f}\colon B\times I \to A(E_{2})$ joining 
$f_{\ast}p^{A}_{1}$ with $p^{A}_{2}$.  Moreover, this homotopy can be lifted to a homotopy 
 joining $f_{\ast}p_{1}^{!}$ with $p_{2}^{!}$  in $Q(E_{2})$ provided that $f$ is homotopic 
to a diffeomorphism of bundles.  The problem of lifting $\omega_{f}$ defines then an obstruction 
to replacing $f$ by a diffeomorphism. This obstruction closely resembles 
the classical Whitehead torsion of a homotopy equivalence of finite CW-complexes. 

Lets return now to the case of a single bundle $p\colon E\to B$ with fiber $M$. 
The map $\lambda \colon Q(E)\to K(R)$ factors though the assembly map so we get a commutative 
diagram 
$$
 \xymatrix{
  &Q(E_{+})\ar[d]^{a} \\
  & A(E) \ar[d] \\
B\ar@(u, u)[uur]^{p^{!}}\ar[ru]^{p^{A}}\ar[r]_{c_{\rho}} & K(R)  \\
}
$$  
If we would have a homotopy equivalence 
$f\colon M\times B \to E$ from the product bundle over $B$ to $p$ then the homotopy 
$\omega_{f}$ would give us a contraction of $p^{A}$  (and thus also a contraction of 
$c_{\rho}$) to a constant map. 
The construction of the smooth torsion takes advantage the fact that under some homological 
conditions the map $c_{\rho}$ is contractible even if we do not have a homotopy equivalence $f$. 
In this case  the contraction of  $c_{\rho}$ 
can be still interpreted as a way of relating the bundle $p$ to the product bundle
on the level of the $K$-theory.
Viewed from this perspective the smooth torsion $\tau_{\rho}(p)$ is  an obstruction 
to the existence of a diffeomorphism between $p$ and the product bundle. 
In effect we can think of it as a linearized Whitehead torsion.
This parallels the way  the classical Reidemeister torsion is interpreted in 
the setting of  CW-complexes.

While the idea of the construction of smooth torsion is simple to explain the technical details 
are rather involved. One of the main problems is that the target of the transfer $p^{!}$ and the 
domain of the map $\lambda$ as described in \cite{DWW} are different and are only linked 
by a zigzag of weak equivalences. 
As a result the torsion of a bundle is not really defined uniquely but rather up to a contractible 
space of choices. This makes this construction of torsion rather inaccessible to direct computations. 
In fact,  at present there are no examples of bundles for which the smooth torsion does 
not vanish, even though such examples abound for Bismut-Lott and Igusa-Klein torsions, and 
intuitively the smooth torsion of Dwyer-Weiss-Williams should be a more delicate invariant. 
The last sentence points out  another problem with the Dwyer-Weiss-Williams construction: at
present there are no known results relating it to the other definitions of torsion of smooth bundles. 

An additional difficulty with the construction of torsion as described above is that it depends on 
existence of the local system of coefficients $\rho$ yielding acyclicity of fibers of $p$. Such 
systems of coefficient are not easy to construct. 

One the goals of this note is to show how these problems can be resolved. 
First, we substantially simplify the construction 
of smooth torsion  using  Waldhausen's manifold approach to the $Q$-construction 
\cite{WalM1}, \cite{WalM2}.  This idea is not entirely original -- in \cite{DWW} the authors sketch it briefly
at the very end of the paper and attribute it to Waldhausen. Our aim, however, is to develop it
 to the extent which would permit us to study properties of the smooth torsion. 
In addition we show that smooth torsion can be defined even if the fibers of $p$ are not acyclic, 
as long as the fundamental group of the base space acts trivially   (or even unipotently -- see \ref{UNIPOTENT}) on the homology groups of the fibers. This demonstrates that  smooth torsion 
exists for a broad class of bundles.

We also aim  to bring the smooth torsion of Dwyer-Weiss-Williams closer to the Bismut-Lott
and the Igusa-Klein constructions. The starting point here is the paper \cite{igusaAx} of Igusa which 
describes a set of axioms for torsion of smooth bundles. Igusa shows that any notion of torsion 
satisfying these axioms must coincide with the Igusa-Klein torsion up to some scalar constants.
In Igusa's setting, torsion is an invariant defined for all smooth unipotent bundles -- 
the condition which as we mentioned above turns out to be  satisfied by the smooth torsion. 
This invariant is supposed to take its values in the cohomology groups $H^{4k}(B, \real)$ 
of the base space of  the bundle. We show that the smooth torsion can be reduced to such a cohomological invariant. 
What remains to be verified is that the cohomology classes we produce satisfy Igusa's axioms. 
This is the goal which the present authors in collaboration with John Klein plan to complete 
in a future paper. 
\end{nn}

\begin{nn}{\bf Organization  of the paper.}
As we have mentioned above our main tool in this paper is the construction of the space 
$Q(X_{+})$ using the language of ``partitions'' given by Waldhausen in \cite{WalM1}, \cite{WalM2}. 
We start by summarizing 
this construction in Section \ref{MANIFOLD}. In \S\ref{ASSEMBLY} 
we describe -- again following Waldhausen --
the assembly map from $Q(X_{+})$ to  Waldhausen's $A$-theory of the space $X$. 
While this map is not  our main interest here, we will use it in Section \ref{TRANSFER} to show
that  a certain map we construct there coincides with the Becker-Gottlieb transfer 
$p^{!}\colon B\to Q(E_{+})$ of a bundle $p\colon E\to B.$
In Section \ref{LINEARIZATION} the assembly map is used again to construct the linearization map 
$\lambda \colon Q(B_{+})\to K(R)$.  
In \S\ref{TORSION} we describe the construction of torsion for bundles with acyclic fibers and 
for unipotent bundles.
Finally, in \S\ref{CH CLASSES} we show how the torsion of a bundle $E\to B$ defines certain 
cohomology classes in $H^{4k}(B; \real),$ for $k>0$. 
\end{nn}


\section{Waldhausen's manifold approach}
\label{MANIFOLD}

\begin{nn}{\bf Partitions.}
Let $I$ denote the closed interval $[0, 1]$. 
For a smooth manifold $X$ with boundary $\partial X$ a partition of $X \times I$ is a triple 
$(M, F, N)$ where $M$ and $N$
are codimension $0$ submanifolds of $X\times I$ such that $M\cup N = X\times I$
and $F=M\cap N$  is a submanifold of codimension $1$. Moreover, 
we assume that $M$ contains $X\times \{0\}$ and is disjoint from $X\times\{1\}$,
and finally, that for some open neighborhood $U$ of $\partial X$ in $X$ the intersection 
of $F$ with $U\times I$ coincides with $U\times \{t\}$ for some $t\in I$.  
While this description reflects the basic properties of partitions we  will make some 
further technical assumptions which will make it easier to work with them:
\begin{itemize}
\item we will assume  that the value of $t$ if fixed once for all partitions (say,  $t=\frac{1}{3}$);
\item we will assume  that $X\times [0, \frac{1}{3}]\subseteq M$  for any partition $(M, N, F)$. 
\end{itemize}
 In the language of \cite{WalM1} partitions satisfying the last two conditions are called lower partitions. 
 
 While the idea is all  components  - $M, N, F$ - of a partition should be smooth
 \begin{floatingfigure}[r]{3.5cm}
\includegraphics[width=3.2cm]{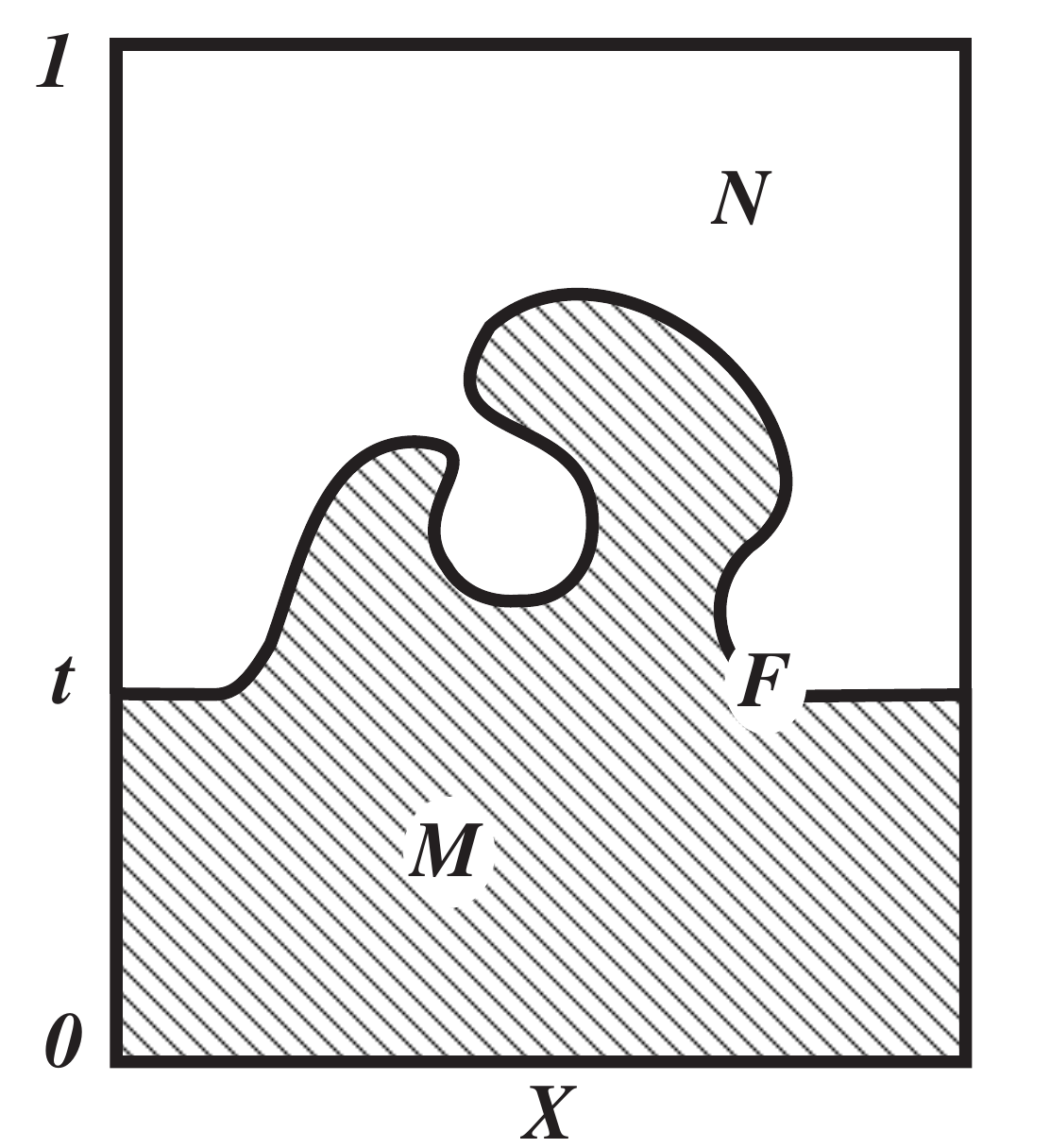}
\end{floatingfigure}
\noindent  this condition
is too rigid for the constructions we will want to perform. Following Waldhausen \cite[Appendix]{WalM1}
we will  assume that these are just topological manifolds but  we also assume  that there is a preferred smooth vector field $s$ on $X\times I$
 which is normal to $F$ in the following
sense.  Given any smooth chart of $X\times I$ containing $x\in F$ there are constants $c, C>0$ 
such that for all $|r|\leq C$ the distance function satisfies the inequality 
$$d(x+rs(x), F) \geq c|r|$$  
This just means that the line passing through $x$ and going in the direction of $s(x)$ stays well 
away from $F$.
\begin{floatingfigure}[l]{0pt} 
\end{floatingfigure}
\parindent 0mm 
The existence of such a vector field $s$ implies that the manifold $F$ admits a smoothing which is 
obtained by sliding points of $F$ along the integral curves of $s$. Moreover, the space of smoothings 
of $F$ which can be obtained in this way is contractible, so we can think of the quadruple 
$(M, N, F, s)$ as describing an essentially unique smooth partition. We again  list some 
additional technical assumptions: 
\begin{itemize}
\item we will assume that on $U\times I$ the vectors of $s$ are the unit vectors pointing 
upward, in the direction of $I$; 
\item we will also assume that for $x\in X\times \{\frac{1}{3}\}$ the component of  $s(x)$ 
tangent to the interval $I$ is a non-zero vector pointing upward. 
\end{itemize}

Clearly a partition $(M, N, F, s)$ is determined by the manifold $M$ and the vector field 
$s$. In order to simplify notation we will write $(M, s)$ instead of $(M, N, F, s)$. 

For a manifold Y we have the notion of a locally trivial family of partitions of $X\times I$ parametrized by
 $Y$. By this we mean a pair $(\bar M, \bar s)$ where $\bar M\subseteq X\times I\times Y$ and 
$s$ is a smooth vector field on $X\times I\times Y$ such that  
\begin{itemize}
\item for each $y\in Y$ the pair $(\bar M\cap X\times I\times \{y\}, \bar s|_{X\times I\times \{y\}})$
is a partition of $X\times I \times \{y\}$
\item for every $y\in Y$ there is an open neighborhood $y\in V\subseteq Y$,
a partition $(M, s)$ of $X\times I$, and a diffeomorphism 
$$\varphi \colon X\times I\times V\to X\times I\times V$$
such that $p_{V}\circ \varphi=p_{V}$ where $p_{V}\colon X\times I\times V\to V$ is the projection 
map, $\varphi$ is the identity map on an appropriate neighborhood of $\partial(X\times I)\times V$, 
$\varphi(M\times V)= \bar M\cap (X\times I\times V)$ and $D\varphi(s)= \bar s$. 
\end{itemize}
If $(\bar M, \bar s)$ is a partition of $X\times I$ parametrized by $Y$ and 
$f\colon Z\ra Y$ is a smooth function then we obtain the induced partition 
$f^{\ast}Y$ parametrized by $Z$.

Denote by $\PP[k](X)$ the set of all partitions parametrized by the $k$-simplex $\Delta^{k}$. 
These sets can be assembled to form a simplicial set $\PP(X)$.  

\end{nn}

\begin{nn}{\bf Stabilization.}
The set $\PP[0](X)$, which 
can be identified with the set of all partitions of $X\times I$, has a partial monoid structure defined 
as follows: given two partitions $(M_{1},s)$ and $(M_{2}, s)$ in $\PP[0](X)$ such that 
$M_{1}\cap M_{2}=X\times [0, \frac{1}{3}]$ we set 
$(M_{1}, s)+(M_{2}, s):=(M_{1}\cup M_{2}, s)$. We can extend this definition to $\PP[k](X)$
for all $k\geq 0$ to obtain a partial simplicial monoid structure on $\PP(X)$. In order to 
define addition for all partitions (and thus define an $H$-space structure on $\PP(X)$) we 
need to introduce stabilization of partitions (called lower stabilization in \cite{WalM1}).  It is a map
of simplicial sets $\sigma\colon \PP(X)\to \PP(X\times J)$ where $J=[0, 1]$. Given a partition 
$(M, s)\in \PP[0](X)$ we set $\sigma(M, s)=(\sigma(M), \sigma(s))$ where 
$\sigma(M)=X\times J\times [0, \frac{1}{3}]\cup M\times [\frac{1}{3}, \frac{2}{3}]$.

\begin{figure}[h]   \includegraphics[width=2.5in]{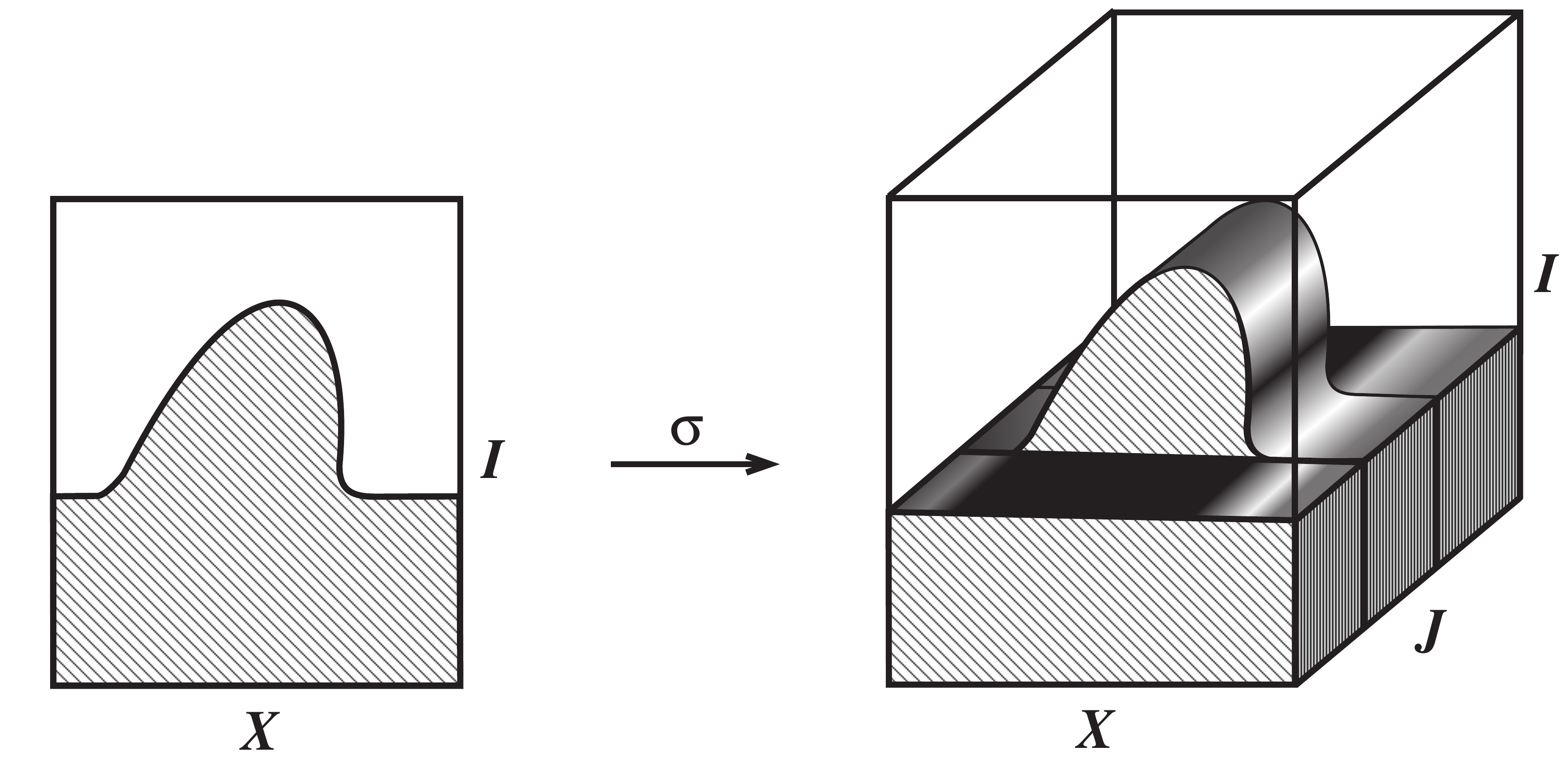} 
  \end{figure}
\noindent In order to define the vector field $\sigma(s)$, fix a smooth vector field $s'$ on 
the interval $J$ such that $s'$ is non-zero at the points $\frac{1}{3}$ and $\frac{2}{3}$
and is zero on some neighborhood of $\partial J$. For $(x, t, t')\in X\times I\times J$ we then set
$$\sigma(s)(x, t, t'):= s(x, t)+ s'(t)$$

We note that the vector field $\sigma(s)$ does not really satisfy all assumptions of our definition 
of a partition since it is not a unit vector field pointing in the direction of $I$ when restricted to a neighborhood of $\partial (X\times J)\times I$. This however can be easily fixed.

In a similar way we can define stabilization maps $\sigma\colon  \PP[k](X)\to \PP[k](X\times J)$
for all $k> 0$ so that we obtain a simplicial  map of  $\sigma \colon \PP(X)\to \PP(X\times J)$.
Notice that given any two partitions it is always possible to slide their 
stabilizations away from each other along $J$  so that the sum is defined. 
In this way the partial monoid structure 
becomes a monoid structure on  $\colim[m] \PP(X\times I^{m})$.
\end{nn}

\begin{nn}{\bf Group completion.}\
The Waldhausen manifold model for $Q(X_{+})$ can be  obtained as a group completion 
of the simplicial monoid $\colim[m]\PP(X\times I^{m})$. In order to describe this group completion 
one can use Thomason's variant of Waldhausen's $S_{\bullet}$-construction (see \cite{Wal1}, p.343). 
Let $\TT[n]\PP[0](X)$ denote the category whose objects are $(n+1)$-tuples  
$\{(M_{i}, s_{i})\}_{i=0}^{n}$ such that $(M_{i}, s_{i})\in \PP[0](X)$, $s_{i}=s_{j}$ for all $i, j$
and that we have inclusions 
$$M_{0}\subseteq M_{1}\subseteq \dots \subseteq M_{n}$$
In $\TT[n]\PP[0](X)$ we have a unique morphism $\{(M_{i}, s_{i})\}\to \{(M'_{i}, s'_{i})\}$ if and only if 
$M'_{0}\cap M_{i}\subseteq M_{0}$ and  $M'_{i}=M_{i}\cup M'_{0}$ for all $i\geq 0$.
Analogously, for any $n, k\geq 0$ we can define a category $\TT[n]\PP[k](X)$
whose objects are increasing sequences of length $n$ in $\PP[k](X)$.
For any fixed $n$ the categories $\TT[n]\PP[k](X)$ assemble to give a simplicial category 
$\TT[n]\PP(X)$. We have functors $d_{j}\colon \TT[n+1]\PP(X) \to 
\TT[n]\PP(X)$ and $s_{j}\colon\TT[n]\PP(X)\to \TT[n+1]\PP(X)$ which are obtained by removing
(or respectively repeating) the $j$-th element of the sequence $\{(M_{i}, s_{i})\}$. In this way
$\TT\PP(X)$ becomes a simplicial object in the category of simplicial categories. 
Let $|\TT\PP(X)|$ denote the space obtained by first taking nerve of each of the categories
$\TT[n]\PP(X)$ and then applying geometric realization to the resulting bisimplicial set. 
Notice that we have a cofibration $|\TT[0]\PP(X)|\to |\TT\PP(X)|$.  One can check that the 
space $|\TT[0]\PP(X)|$ is contractible,
and so we get a weak equivalence 
$$|\TT\PP(X)|\simeq |\TT\PP(X)|/|\TT[0]\PP(X)|$$
By abuse of notation from now on we will denote by $|\TT\PP(X)|$ the quotient space 
on the right. The advantage of this modification is that $|\TT\PP(X)|$ has now a canonical 
choice of a basepoint. 

Since everything we did so far behaves well with respect to the stabilization maps we 
obtain the induced maps of spaces 
$$\sigma\colon  |\TT\PP(X\times I^{m})|\to |\TT\PP(X\times I^{m+1})|$$
Passing to the homotopy colimit we get
\begin{theorem}[Waldhausen\cite{WalM2}]
There is a weak equivalence
$$\Omega\, \hocolim[m]|\TT\PP(X\times I^{m})| \simeq Q(X_{+})$$
\end{theorem}
In view of this result from now on by $Q(X_{+})$ we will understand the space on the left
hand side of the above equivalence.
\end{nn}

\begin{remark}
\label{TT1}
The following observation will be useful for constructing maps into  $Q(X_{+})$.
Notice that we have a map 
$$|\TT[1]\PP(X)|\times \Delta^{1}\to |\TT\PP(X)|$$
After stabilizing the right hand side and taking the adjoint 
we obtain a map $|\TT[1]\PP(X)|\to Q(X_{+})$. Thus any map into the nerve of the 
category $\TT[1]\PP(X)$ naturally yields a map into $Q(X_{+})$.
\end{remark}

\section{The assembly map}
\label{ASSEMBLY}
Waldhausen's motivation for constructing the space $Q(X_{+})$  in the way sketched in the last 
section was to relate this space to $A(X)$ -- the $A$-theory of the space $X$. Since we will
need to use this relationship later in this paper, we now  describe a map $a\colon Q(X_{+})\to A(X)$
which we will call the assembly map. 

The simplest way of constructing the space $A(X)$ is to start with the category $\Rfd(X)$
whose objects are homotopy finitely dominated retractive spaces over $X$ and whose 
morphisms are maps of retractive spaces. The category $\Rfd(X)$ is a Waldhausen category 
in the sense of \cite[ Definition 1.2]{Wal1} with cofibrations given by Serre cofibrations and weak equivalences 
defined as weak homotopy equivalences. It follows that we can turn it into a simplicial category 
$\TT\Rfd(X)$ using again the Thomason's variant of Waldhausen's $\SS$-construction. We  set
$$A(X):=\Omega (|\TT\Rfd(X)|/|\TT[0]\Rfd(X)|)$$

In order to obtain a direct map from our model of $Q(X_{+})$ this construction needs to 
be modified somewhat. 
First, for $k\geq 0$ we can construct a category $\Rfd_{k}(X)$ whose objects are locally homotopy 
trivial families of retractive spaces over $X$ parametrized by the simplex $\Delta^{k}$.
These categories taken together form a simplicial category $\Rfd_{\bullet}(X)$. 
Analogously as we did in the case of categories of partitions we can define stabilization functors
$$\sigma\colon \Rfd_{\bullet}(X)\to \Rfd_{\bullet}(X\times I)$$
in the following way: if $Y$ is a retractive space over $X$ then 
$$\sigma(Y) := Y\times [\frac{1}{3}, \frac{2}{3}]\cup_{X\times [\frac{1}{3}, \frac{2}{3}]} X\times I$$

Applying the $\TT$-construction to $\Rfd_{\bullet}(X)$ we get a bisimplicial category 
$\TT\Rfd_{\bullet}(X)$. Define 
$$A_{p}(X):= \hocolim[m] |\TT\Rfd(X\times I^{m})|/|\TT[0]\Rfd(X\times I^{m})|$$

Notice that if $X$ is a smooth manifold and if $(M, s)$ is a partition of $X\times I$ then 
$a(M):=M\cap (X\times [\frac{1}{3}, 1])$ is in an obvious way a retractive space over $X$. 
The assignment $(M, s)\mapsto a(M)$ extends to a functor of simplicial categories 
$a\colon \PP(X)\to \Rfd_{\bullet}(X)$ which commutes with the stabilization functors. 
This induces a map 
$$a\colon Q(X_{+})\to A_{p}(X)$$

Since the category $\Rfd(X)$ is isomorphic to $\Rfd_{0}(X)$ we have a  functor 
$i\colon \Rfd(X)\to \Rfd_{\bullet}(X)$ which induces a map $i\colon A(X)\to A_{p}(X)$. We have

\begin{theorem}[Waldhausen,  {\cite[Lemma 5.4]{WalM1}}]
\label{A-AP}
The map $i\colon A(X)\to A_{p}(X)$ is a homotopy equivalence.
\end{theorem} 
 
 While this result says that we are not changing much by replacing $A(X)$ with $A_{p}(X)$,
 it will be convenient
 to have an assembly map whose codomain is $A(X)$. In order to get such a map 
define  $\tilde Q(X_{+})$ to be  the homotopy pullback of the diagram 
\begin{equation}
\label{QTILDE}
\xymatrix{A(X)\ar[r] & A_{p}(X)   & Q(X_{+})\ar[l]_{a}}
\end{equation}
In view of Theorem \ref{A-AP} we have $\tilde Q(X_{+})\simeq Q(X_{+})$, and 
$\tilde Q(X_{+})$ comes equipped with a map
$\tilde a \colon \tilde Q(X_{+})\to A(X)$.

\section{Transfer}
\label{TRANSFER}
Going back to the diagram on page \pageref{MAIN_DIAG} we see that if we want to describe 
the smooth torsion of a bundle $p\colon E\to B$ we need to construct the transfer map
$p^{!}\colon B \to Q(E_{+})$ (or rather $p^{!}\colon B \to \tilde Q(E_{+})$) and the linearization 
map $\lambda \colon  \tilde Q(E_{+})\to K(R)$.  We deal with the transfer in this section and with the linearization in the next one.

Let $p\colon E\to B$ be a smooth bundle of manifolds with $B$ and $E$ compact. Denote by $T^{v}E$
the subbundle of the tangent bundle $TE$ consisting of vectors tangent to the fibers 
of $p$. By choosing a Riemannian metric on $TE$ we can identify the bundle $p^{\ast}TB$ with the 
subbundle of $TE$ which is the orthogonal complement of $T^{v}E$. Using this identification
given the map $p\times {\id}\colon E\times I\to B\times I$ we can consider the bundle 
$(p\times \id)^{\ast}T(B\times I)$ as a subbundle of $T(E\times I)$.
As a consequence any section $s\colon B\times I\to T(B\times I)$ will define a section 
$(p\times \id)^{\ast}s$ of the bundle $T(E\times I)$.   

Assume for a moment  that fibers of $p$ are closed manifolds. In this case given a partition 
$(M, s)\in \PP[0](B)$ the pair $((p\times\id)^{-1}M, (p\times \id)^{\ast}s)$ defines a partition of $E\times I$, 
so we get a map
of sets $Q(p^{!})\colon\PP[0](B)\to \PP[0](E)$. Since this map preserves the partial ordering 
of partitions  we in fact obtain a functor $\TT[0]\PP[0](B)\to \TT[0]\PP[0](E)$.  
In the same way we can define functors $\TT[n]\PP[k](B)\to \TT[n]\PP[k](E)$  for all 
$k, n\geq 0$ so that they induce a map
$$Q(p^{!})\colon |\TT\PP(B)|\to |\TT\PP(E)|$$
Since $Q(p^{!})$ is compatible with stabilization we get a map $Q(p^{!})\colon Q(B_{+})\to Q(E_{+})$.

If fibers of the bundle $p$ are manifolds with boundary we need to modify the above construction 
slightly so that for a partition $(M, s)\in \PP[0](B)$ the element $Q(p^{!})(M, s)$ behaves nicely near 
$\partial E \times I$.  This can be done as follows. Let $\partial^{v}E$ be the subspace
of $E$ consisting of all boundary points of the fibers of $p$. If $F$ is a fiber of $p$ then 
$p|_{\partial^{v}E}\colon \partial^{v} E\to B$ is a subbundle of $p$ with fiber $\partial F$. 

 For $b\in B$ let $F_{b}$
denote the fiber of $p$ over $b$. We can find an open neighborhood $U\subseteq E$ in such way that 
for all $b\in B$ the intersection $U\cap F_{b}$ is a collar neighborhood of 
$\partial F_{b}$ in $F_{b}$ \cite[p.590]{Bec}. 
For a partition $(M, s)\in \PP[0](M)$ let $Q(p^{!})(M)\subseteq E\times I$ be given by
$$Q(p^{!})(M):= U\times [0, \frac{1}{3}] \cup ((p\times \id)^{-1}(M)\cap ((E\setminus U)\times I)) $$
We need then to modify the vector field $(p\times \id)^{\ast}s$ so that it is normal to $Q(p^{!})(M)$. 
This can be done in a way similar to the one we used to define stabilization of partitions. 

In a similar way given a bundle $p\colon E\to B$ we can construct maps 
$A(p^{!})\colon A(B)\to A(E)$ and $A_{p}(p^{!})\colon A_{p}(B)\to A_{p}(E)$. Each of these maps 
is induced by a functor of categories of retractive spaces which assigns to a retractive space 
over $B$ its pullback along $p$ (in order to make this compatible with the construction of $Q(p^{!})$ 
we need to modify these pullback slightly in a neighborhood of the boundaries of fibers of 
$p$). These three maps induce in turn a map of homotopy pullbacks
$$\tilde Q(p^{!})\colon \tilde Q(B_{+})\to \tilde Q(E_{+})$$

The map $p^{!}\colon B\to \tilde Q(E_{+})$ will be obtained as the composition of $\tilde Q(p^{!})$
with a coaugmentation map $\eta\colon B\to \tilde Q(B_{+})$ which we describe below. For simplicity 
we will assume first that $B$ is a closed manifold.

Let $\simp$ denote the simplicial set of singular simplices of $B$. It will be convenient to consider it 
as a simplicial category with identity morphisms only on each simplicial level. 
We have a weak equivalence $B\simeq |\simp|$, so it will suffice to construct a map 
$|\simp|\to \tilde Q(B_{+})$. Recall that  $\tilde Q(B_{+})$ was defined as the homotopy 
pullback of the diagram  (\ref{QTILDE}). 
Notice also  that we have a commutative diagram 
$$
\xymatrix{
|\TT[1]\Rfd(B)| \ar[r]\ar[d] & |\TT[1]\Rfd_{\bullet}(B)|\ar[d] & |\TT[1]\PP(B)| \ar[l]\ar[d] \\
A(B) \ar[r] & A_{p}(B)  & Q(B_{+}) \ar[l] \\
}
$$
where the vertical maps are obtained as in  Remark \ref{TT1}. It will then suffice to define a map 
from $|\simp|$ to the homotopy limit of the upper row of this diagram. This, in turn, can be 
accomplished by specifying functors $ \simp \to \TT[1]\Rfd(B)$ and $ \simp \to \TT[1]\PP(B)$ 
and a zigzag of natural transformations joining these functors in $\TT[1]\Rfd_{\bullet}(B)$. 
A minor inconvenience here is the fact that $\TT[1]\PP(B)$,  $\TT[1]\Rfd_{\bullet}(B)$  and 
$\simp$ are simplicial categories while  $\Rfd(B)$ is not, but we can think about $\Rfd(B)$
as of a constant simplicial object in the category of small categories. 

We will start then with a diagram of categories 
$$\xymatrix{
\TT[1]\Rfd(B)\ar[r]^{i} & \TT[1] \Rfd_{\bullet}(B) & \TT[1]\PP(B) \ar[l]_{a} \\
}$$
and we will extend it  to a diagram 
$$\xymatrix{
 & \simp \ar[dl]_{\eta_{\RR}} \ar[d]^{\eta_{p\RR}}\ar[dr]^{\eta_{\PP[]}} \\
\TT[1]\Rfd(B)\ar[r]^{i} & \TT[1] \Rfd_{\bullet}(B) & \TT[1]\PP(B) \ar[l]_{a} \\
}$$
such that the two triangles of functors commute up to natural transformations.

The functor $\eta_{\RR}\colon \simp\to \TT[1]\Rfd(B)$ is defined as follows: 
given a singular simplex $\sigma\colon \Delta^{k}\to B$ consider the  retractive space
$B\sqcup \Delta^{k}$ over $B$.  We set $\eta_{\RR}(\sigma)$ to be the cofibration 
of retractive spaces $B \hra B \sqcup \Delta^{k}$. 

In order to define the functor $\eta_{\PP[]}\colon \simp \to \TT[1]\PP(B)$
fix a Riemannian metric on the tangent bundle $TB$. Choose $\epsilon >0$ such that
the closed disc bundle $TB_{\epsilon}$ consisting of vectors of $TM$ of length $\leq \epsilon$
has the property that the exponential map ${\rm exp}\colon TB_{\epsilon} \to M$ restricted 
to each fiber of $TB_{\epsilon}$ is a diffeomorphism onto its image. 
Given a singular simplex $\sigma\colon \Delta^{k}\to B$ consider the induced disc bundle
$\sigma^{\ast}TB_{\epsilon}$ over $\Delta^{k}$. The exponential map gives a map of bundles
$$ 
\xymatrix{
\sigma^{\ast}TB_{\epsilon} \ar[rr]^{\rm exp} \ar[rd] & & B\times \Delta^{k} \ar[ld] \\
& \Delta^{k} & \\
}$$
which is a fiberwise embedding. Fix numbers $a, b$ such that $\frac{1}{3}< a < b <1$. 
We have a fiberwise embedding of fiber bundles over $\Delta^{k}$:
$$ 
\xymatrix{
(\sigma^{\ast}TB_{\epsilon}\times [a, b]) \cup (B\times [0, \frac{1}{3}] \times \Delta^{k})\ar[rr] \ar[rd] 
& & B\times I \times \Delta^{k} \ar[ld] \\
& \Delta^{k} & \\
}$$
This fiberwise embedding defines a family of partitions parametrized by $\Delta^{k}$. 
We set $\eta_{\PP[]}(\sigma)$ to be the inclusion of families of partitions
$$B\times [0, \frac{1}{3}]\times \Delta^{k} \hra (\sigma^{\ast}TB_{\epsilon}\times [a, b]) \cup (B\times [0, \frac{1}{3}]\times \Delta^{k})$$

Next, we need to define the functor $\eta_{p\RR}\colon \simp \to \TT[1]\Rfd_{\bullet}(B)$.
For a singular simplex $\sigma\colon \Delta^{k}\to B$  we have the map 
$$\id_{\Delta^{k}} \sqcup {\rm pr}_{\Delta^{k}} \colon \Delta^{k} \sqcup B\times \Delta^{k}\to \Delta^{k}$$ 
For each $t\in \Delta^{k}$ the fiber of this map over $t$ is in a natural way a retractive space over 
$B$, so we can think of $ \Delta^{k} \sqcup B\times \Delta^{k}$ as of a family of retractive spaces
over $B$ parametrized by $\Delta^{k}$. The functor 
$\eta_{p\RR}$ is given by the assignment
$$\eta_{p\RR}(\sigma) := (B\times \Delta^{k}\hra \Delta^{k} \sqcup B\times \Delta^{k} )$$ 

In order to describe a natural transformation from $\eta_{p\RR}$ to $i\eta_{\RR}$  notice that 
for $\sigma\colon \Delta^{k}\to B$ we obtain $i\eta_{\RR}(\sigma)$ 
by taking the retractive space $\eta_{\RR}(\sigma)= \Delta^{k}\sqcup B$ and multiplying it by 
$\Delta^{k}$ which makes it into a retractive space over $B$ parametrized by $\Delta^{k}$. 
The natural transformation is then defined by the maps
$$\eta_{p\RR}(\sigma)= \Delta^{k}\sqcup  B\times \Delta^{k} \to 
(\Delta^{k}\times\Delta^{k})\sqcup (B\times \Delta^{k})=(\Delta^{k}\sqcup B)\times \Delta^{k}
=i\eta_{\RR}(\sigma)$$
which restrict to the identity map on $B\times \Delta^{k}$ and which send $x\in \Delta^{k}$ to 
$(x, x)\in \Delta^{k}\times\Delta^{k}$.  The natural transformation from $a\eta_{\PP[]}$ to 
$\eta_{p\RR}$ is easy to define.

If $B$ is a manifold with a boundary we need to modify this construction somewhat so that 
the values of the functor $\eta_{\PP}$ are still partitions. This can be done by choosing an 
open collar neighborhood $U$ of the boundary of $B$. On $B\setminus U$ the map 
$\eta$ can be now defined exactly as before. We then extend it to $B$ by composing it 
with the retraction $B\to B\setminus U$. 

\begin{proposition}
If $p\colon E\to B$ is a smooth fibration then the  map 
$$p^{!}:=\tilde Q(p^{!})\circ \eta \colon B \to Q(E_{+})$$
is the Becker-Gottlieb transfer of $p$.  
\end{proposition}

\begin{proof}
The composition $\tilde a \circ p^{!}$ is homotopic to $\ehom(p)$ - 
the homotopy Euler characteristic of the bundle $p$ as defined 
in \cite{DJ}. By \cite[Thm. 3.12]{DJ} we obtain that ${\rm tr}\circ \ehom(p)$  is 
the Becker-Gottlieb transfer where ${\rm tr}\colon A(E)\to \tilde Q(E_{+})$ is the 
Waldhausen's trace map \cite{WalM2}. By \cite{WalM2} we have   
${\rm tr} \circ \tilde a\sim \id_{\tilde Q(E_{+})}$, so 
$${\rm tr}\circ \ehom(p)\sim {\rm tr}\circ \tilde a \circ p^{!} \sim p^{!}$$
\end{proof}

\section{Linearization}
\label{LINEARIZATION}

For a ring  $R$ let  $\Ch(R)$ denote the category of finitely homotopy dominated chain complexes of projective $R$-modules.  The category $Ch(R)$ can be equipped with a Waldhausen model 
category structure with degreewise monomorphisms as cofibrations and quasi isomorphisms 
as weak equivalences. Applying Waldhausen's $\SS$-construction we obtain a simplicial 
category $\SS Ch(R)$. The associated space $K(R)=\Omega (|\SS Ch(R)|)$
is homotopy equivalent to the infinite loop space underlying the $K$-theory spectrum 
of the ring $R$.

Let $X$ be a space and let $\rho \colon {\mathcal M} \to X$ be a locally constant sheaf  
of finitely generated projective $R$-modules.  
As in \cite[p.40]{DWW} we notice that we have a functor 
$\lambda^{\RR}_{\rho}\colon \Rfd(X)\to Ch(R)$
which assigns to every retractive space $Y\in \Rfd(X)$ the  relative singular chain complex 
of  $C_{\ast}(Y, X,  \rho)$ with coefficients in $\rho$.  This functor induces a map 
$\lambda^{\RR}_{\rho}\colon A(X)\to K(R)$. Recall that for a smooth manifold $X$ we 
constructed the assembly map $\tilde a \colon \tilde Q(X_{+})\to A(X)$. 
By the  linearization map $\lambda_{\rho}\colon \tilde Q(X)\to K(R)$ 
we will understand the composition  $\lambda^{\RR}_{\rho}\circ \tilde a$.

\section{Smooth torsion}
\label{TORSION}
We are now in position to define smooth Reidemeister torsion of a 
bundle of manifolds. We will do it under two different sets of assumptions, 
one replicating the conditions of \cite{DWW}, and the other conforming to the 
axiomatic setup of Igusa \cite{igusaAx}. We note, however, that the idea underlying both 
constructions is essentially the same: given a bundle $p\colon E\to B$ and 
a sheaf of $R$-modules $\rho\colon {\mathcal M} \to E$ we have constructed maps
$p^{!}\colon B \to \tilde Q(E_{+})$ and $\lambda_{\rho}\colon \tilde Q(E_{+})\to K(R)$. 
Consider the composition 
$$\lambda_{\rho}\circ p^{!} \colon |\simp|\to K(R)$$
Under certain conditions on the bundle $p$ and the sheaf $\rho$ this map is 
homotopic to a constant map via a preferred homotopy. As a consequence 
we obtain a lift of $p^{!}$ to a map 
$\tsm(p)\colon B \to \hofib(\tilde Q(E_{+})\to K(R))$. 
This lift is the smooth torsion of the bundle $p$.

\begin{nn}{\bf Acyclicity.} 
Assume that we are given a sheaf $\rho$ such that for any $b\in B$ we 
have $H_{\ast}(F_{b}; \rho)=0$ where $F_{b}$ is the fiber of $p$ over $b$. 
Notice that the map $\lambda_{\rho}\circ p^{!}$ comes from a functor 
$\simp\to Ch(R)$ which assigns to each simplex $\sigma\colon \Delta^{k}\to B$
the relative chain complex  $C_{\ast}(\sigma^{\ast} E \sqcup E, E; \rho)$ where 
$\sigma^{\ast} E$ denotes the pullback of the diagram 
$$
\xymatrix{
\Delta^{k} \ar[r]^{\sigma} & B & E \ar[l]_{p} \\
}
$$
Vanishing of homology groups of the fibers of $p$ implies that all these chain complexes 
are acyclic, and thus the maps $C_{\ast}(\sigma^{\ast} E \sqcup E, E; \rho)\to 0$
are all quasi isomorphisms and define a natural transformation to the functor which sends
all simplices of $\simp$ to the zero chain complex. On the level of spaces 
this natural transformation defines a homotopy $\omega_{\rho}\colon |\simp| \times I \to K(R)$. 
Denote by $\wh_{\rho}(E)$ the homotopy fiber of the linearization map $\lambda_{\rho}$
taken over the basepoint of $K(R)$ represented by the zero chain complex. 
The smooth torsion of the bundle $p$ is the map $\tsm(p)\colon |\simp|\to \wh_{\rho}(E)$
determined by the transfer $p^{!}$ together with the homotopy $\omega_{\rho}$.  

\end{nn}

\begin{nn}{\bf Unipotent bundles.}
\label{UNIPOTENT}
Let $\fd$ be a field and let $p\colon E\to B$ be a bundle such that $B$ is a connected manifold
with a basepoint $b_{0}\in B$. Assume that the fundamental group $\pi_{1}(B, b_{0})$
acts trivially on the homology $H_{\ast}(F_{b_{0}}; \fd)$ of the fiber over $b_{0}$. 

Consider the map $\lambda_{\rho}\colon \tilde Q(E_{+})\to K(\fd)$ associated with the constant sheaf
of 1-dimensional vector spaces over $\fd$. In this case the composition  
$\lambda_{\rho}\circ p^{!}$ assigns to a simplex $\sigma$ the chain complex 
$C_{\ast}(\sigma^{\ast}E\sqcup E, E; \fd)$. We will construct a sequence of 
homotopies joining this map with the constant function which maps the whole space
$|\simp|$ to  $H_{\ast}(F_{b_{0}}; \fd)$, which we will consider as a chain complex with
the trivial differentials. 

Let $wCh({\fd})$ denote the subcategory of $Ch({\fd})$ with the same 
objects as $Ch({\fd})$, but with quasi-isomorphisms as morphisms. We have the canonical 
map 
$$k\colon |wCh({\fd})|\to K({\fd})$$ 
Notice that the map $\lambda_{\rho}\circ p^{!}$
admits a factorization 
$$
\xymatrix{
|\simp| \ar[rr]^{\lambda_{\rho}\circ p^{!}} \ar[dr]& & K(\fd)\\
&|wCh(\fd)| \ar[ur]_{k} & }
$$

We have

\begin{lemma}[{\cite[Prop. 6.6]{DWW}}]
Let $H\colon |wCh(\fd)| \to K(\fd)$ be the map which assigns to each chain complex 
$C$ its homology complex $H_{\ast}(C)$. There is a preferred homotopy $k\simeq H$
\end{lemma}

\begin{proof}
For a chain complex 
$$
\xymatrix{
C = (\dots \ar[r]  & C_{2}\ar[r]^{\partial_{2}} & C_{1}\ar[r]^{\partial_{1}} & C_{0}\ar[r] & 0 )
}
$$ 
Let $P_{q}C$ denote the complex such that $(P_{q}C)_{i} =0$ for $i>q+1$, 
$(P_{q}C)_{q+1}= \partial(C_{q+1})$, and $(P_{q}C)_{i}=C_{i}$ for $i\leq q$. 
Let $Q_{q}C$ be the kernel of the map $P_{q}C\to P_{q+1}C$. We obtain cofibration
sequences 
$$Q_{q}C\to P_{q}C\to P_{k-1}C$$
functorial in $C$. 
Notice that 
the complex $Q_{q}C$ is quasi-isomorphic to its homology complex $H_{\ast}(Q_{q}C)$
and this last complex has only one non-zero module $H_{q}(C)$ in the degree $q$. 
By Waldhausen's additivity theorem we obtain that the map $P_{q}\colon |wCh(\fd)|\to K(\fd)$
which assigns to $C\in wCh(\fd)$ the chain complex $P_{q}C$ is homotopic to the map
$H_{q}\colon |wCh(\fd)|\to K(\fd)$ sending $C$ to $P_{q-1}C \oplus H_{\ast}(Q_{q}C)$.
Iterating this argument we see that for each $q$ the map $P_{q}$ is homotopic to the map 
which assigns to a complex $C$ the chain complex $\bigoplus_{i=0}^{q}H_{i}(C)$.  
Since $C= \lim_{q}P_{q}C$ we obtain the statement of the lemma.
\end{proof}

As a consequence of the lemma we get a homotopy between 
$\lambda_{\rho}\circ p^{!}$ and the map which assigns to a simplex $\sigma$
the homology chain complex  $H_{\ast}(\sigma^{\ast}E\sqcup E, E; \fd)$. 
Since this last chain complex is  isomorphic 
to the  chain complex $H_{\ast}(\sigma^{\ast}E ; \fd)$ 
 we obtain a homotopy from $\lambda_{\rho}\circ p^{!}$ to the map represented by a functor 
 $v\colon \simp \to Ch(R)$ given by  $v(\sigma)=H_{\ast}(\sigma^{\ast}E; \fd)$.

Next, let $v_{0}\colon \simp \to Ch(R)$ denote the 
functor this assigns to each singular simplex $\sigma$ the complex 
$H_{\ast}(F_{\sigma(0)}; \mathbb{F})$
where $F_{\sigma(0)}$ is the fiber of the
bundle $p$ taken over the zero vertex of $\sigma$.  The isomorphisms 
$H_{\ast}(F_{\sigma(0)}; \fd)\to H_{\ast}(\sigma^{\ast}E; \fd)$ form a natural transformation 
of functors $v$ and $v_{0}$. Finally, given any point $b\in B$ choose a path 
joining this point to the basepoint $b_{0}$. Lifting it to the space $E$ we can produce 
a homotopy equivalence $F_{b}\to F_{b_{0}}$. The map which it induces on 
the homology groups will not depend on the choice of the path by our assumption 
the $\pi_{1}(B)$ acts trivially on the homology of the fibers.
The maps $H_{\ast}(F_{\sigma}(0); \fd)\to H_{\ast}(F_{b_{0}}; \fd)$ yield the natural 
transformation from $v_{0}$ to the constant functor. 

On the level of spaces the natural transformations of functors we described above define 
a homotopy from the map $\lambda_{\rho}\circ p^{!}\colon |\simp|\to K(\fd)$ 
to the constant map which maps $|\simp|$ to the point of $K(\fd)$ represented by 
the chain complex $H_{\ast}(F_{b_{0}}; \fd)$. This homotopy taken together 
with the transfer map $p^{!}\colon |\simp|\to \tilde Q(E_{+})$ defines a map
$\tilde \tau_{\fd}(p)\colon |\simp|\to \hofib(\tilde Q(E_{+})\to K(\fd))_{H_{\ast}(F_{b_{0}}, \fd)}$. 
We will call this element the unreduced Reidemeister torsion of the bundle $p$. 
The obvious inconvenience of this definition is that changing a basepoint in $B$
changes the target of the map $\tilde \tau_{\fd}(p)$. This can be fixed by shifting this 
map so it takes values in the space $\wh_{\fd}(E):=\hofib(\tilde Q(E_{+}\to K(\fd))_{0}$
-- the homotopy fiber taken over the zero chain complex. Since both 
$Q(E_{+})$ and $K(\fd)$ are infinite loop spaces this shift can be accomplished by
subtracting the element $p^{!}(b_{0})$ from the map $p^{!}$ and subtracting 
$H_{\ast}(F_{b_{0}}; \fd)$ from the contracting homotopy $|\simp|\times I\to K(\fd)$. 
One could make it more explicit by constructing models for inverses of elements 
in $\tilde Q(E_{+})$ and $K(\fd)$. This is not hard to do. 
The new map $\tau_{\fd}(p)\colon |\simp|\to \wh_{\fd}(E)$ is the (reduced) torsion 
of $p$. 

The above construction can be also carried out under more general conditions which  conform
to the setting of  \cite{igusaAx}.  
We will say that a smooth bundle $p\colon E\to B$ is unipotent if the homology groups 
$H_{\ast}(F_{b_{0}}; \fd)$ admit a filtration
$$0=V_{0}(F_{b_{0}})\subseteq V_{1}(F_{b_{0}})\subseteq \dots V_{k}(F_{b_{0}})=H_{\ast}(F_{b_{0}}; \fd)$$
such that $\pi_{1}(B)$ acts trivially on the quotients $V_{i}/V_{i-1}$. In this case 
consider the functor $v_{0}\colon \simp \to Ch(\fd)$ defined above. Waldhausen's 
additivity theorem implies that the map $v\colon |\simp|\to K(R)$ is canonically 
homotopic to the map which assigns to a simplex $\sigma$ the direct sum 
$\bigoplus V_{i}(F_{\sigma(0)})/V_{i-1}(F_{\sigma(0)})$. Triviality of the action of 
the fundamental group of $B$ on the quotients $V_{i}(F_{\sigma(0)})/V_{i-1}(F_{\sigma(0)})$
implies that we can construct a map $\tau_{\fd}\colon |\simp|\to \wh_{\fd}(E)$ similarly 
as before.

\section{Characteristic classes}
\label{CH CLASSES}

As we mentioned at the beginning of this paper the torsion invariants of smooth bundles 
constructed by Igusa-Klein and Bismut-Lott are constructed as certain cohomology classes 
associated to the bundle. More precisely, for a bundle $p\colon E\to B$ its torsion in 
both of these settings is an element of $\bigoplus_{k>0}H^{4k}(B; {\mathbb R})$. 
Our final goal in this note is to show that the construction of torsion described above 
also gives rise to an element of $\bigoplus_{k>0}H^{4k}(B; {\mathbb R})$ which
brings it on a common ground with the other notions of torsion. 

Let then $p\colon E\to B$ be a bundle with a unipotent action of $\pi_{1}(B)$ on the homology 
of the fiber $H_{\ast}(F_{b_{0}}; {\mathbb R})$, so that the torsion 
$\tau_{\mathbb R}(p)\colon |\simp| \to \wh_{\mathbb R}(E)$ is defined. 
Consider an embedding $i\colon E\to D^{N}$ where $D^{N}$ is a closed disc in ${\mathbb R}^{M}$
for some large $M>0$. Let $NE$ be a closed tubular neighborhood of $E$ in $D^{N}$. Considering 
$NE$ as a disc bundle over $E$ we obtain a transfer map $\tilde Q(E_{+})\to \tilde Q(NE_{+})$. 
Also, since the construction of $\tilde Q$ is functorial with respect to codimension 0 embeddings 
of manifolds  the inclusion $NE\hra D^{M}$ induces a map $\tilde Q(NE)\to \tilde Q(D^{M})$. 
Composing it with the transfer of the bundle $NE\to E$ we obtain a map 
$\tilde Q(E_{+})\to \tilde Q(D^{M}_{+})$. 
Consider the diagram 
$$
\xymatrix{
 & \wh_{\mathbb R}(E)\ar[d]  \ar@{-->}[r]& \wh_{\mathbb R}(D^{M})\ar[d]\\
B\ar[ru]^{\tau_{\mathbb R}(p)} \ar[r]^{p^{!}} & \tilde Q(E_{+})\ar[d]_{\lambda_{\mathbb R}} \ar[r]
& \tilde Q(D^{M}_{+})\ar[ld]^{\lambda_{\mathbb R}} \\
 & K({\mathbb R}) & \\
}
$$
One can check that the lower triangle commutes up to a preferred choice of homotopy, 
so that we obtain a map of homotopy fibers $\wh_{\mathbb R}(E)\to \wh_{\mathbb R}(D^{M})$. 

Now, consider the fibration sequence 
$$\Omega K(\mathbb R)\to \wh_{\mathbb R}(D^{M})\to \tilde Q(D^{M}_{+})\to K({\mathbb R})$$ 
Since this is a fibration of infinite loop spaces after applying the rationalization functor 
we obtain a new fibration sequence.
$$\Omega K(\mathbb R)_{\mathbb  Q}\to \wh_{\mathbb R}(D^{M})_{\mathbb Q}\to 
\tilde Q(D^{M}_{+})_{\mathbb Q}\to K({\mathbb R})_{\mathbb Q}$$
 
Since the  homotopy groups of $ \pi_{i}\tilde Q(D^{M}_{+})\cong \pi_{i}Q(S^{0})$ are torsion 
for $i>0$ the space $\tilde Q(D^{M}_{+})$ is homotopically discrete. It follows that every 
connected component of $\wh_{\mathbb R}(D^{M})_{\mathbb Q}$ is weakly equivalent 
to the space $\Omega K({\mathbb  R})_{\mathbb Q}$.  On the other hand by \cite{Borel} 
we have a weak equivalence 
$$K(\mathbb R)_{\mathbb Q}\simeq {\mathbb Z}\times \prod_{k>0} K({\mathbb R}, 4k+1)$$
Let $\wh_{\mathbb  R}(D^{M})^{B}_{\mathbb Q}$ denote the connected component of 
$\wh_{\mathbb R}(D^{M})_{\mathbb Q}$ which is a target of the map 
$B\to \wh_{\mathbb R}(D^{M})_{\mathbb Q}$. By the observation above we have 
$\wh_{\mathbb R}(D^{M})^{B}_{\mathbb Q}\simeq \prod_{k>0}K({\mathbb R}, 4k)$, 
and so the homotopy class of the map $B\to \wh_{\mathbb R}(D^{M})_{\mathbb Q}$
determines an element in $\bigoplus_{k>0}H^{4k}(B; {\mathbb R})$.

\end{nn}

\bibliographystyle{plain}

\begin{thebibliography}{10}

\bibitem{Bec}
James~C. Becker and Reinhard~E. Schultz.
\newblock Axioms for bundle transfers and traces.
\newblock {\em Math. Z.}, 227(4):583--605, 1998.

\bibitem{Bis-Lott}
Jean-Michel Bismut and John Lott.
\newblock Flat vector bundles, direct images and higher real analytic torsion.
\newblock {\em J. Amer. Math. Soc.}, 8(2):291--363, 1995.

\bibitem{Borel}
Armand Borel.
\newblock Stable real cohomology of arithmetic groups.
\newblock {\em Ann. Sci. \'Ecole Norm. Sup. (4)}, 7:235--272 (1975), 1974.

\bibitem{DJ}
Wojciech Dorabiala and Mark~W. Johnson.
\newblock Factoring the {B}ecker-{G}ottlieb transfer through the trace map.
\newblock {\em {\rm Math arXiv} \tt math/0601620}, 2006.

\bibitem{DWW}
W.~Dwyer, M.~Weiss, and B.~Williams.
\newblock A parametrized index theorem for the algebraic {$K$}-theory {E}uler
  class.
\newblock {\em Acta Math.}, 190(1):1--104, 2003.

\bibitem{Igusa}
Kiyoshi Igusa.
\newblock {\em Higher {F}ranz-{R}eidemeister torsion}, volume~31 of {\em AMS/IP
  Studies in Advanced Mathematics}.
\newblock American Mathematical Society, Providence, RI, 2002.

\bibitem{igusaAx}
Kiyoshi Igusa.
\newblock Axioms for higher torsion invariants of smooth bundles.
\newblock {\em {\rm To appear in the} Journal of Topology, {\rm Math arXiv} \tt
  math/0503250}, 2005.

\bibitem{Klein}
John~R. Klein.
\newblock Higher {R}eidemeister torsion and parametrized {M}orse theory.
\newblock In {\em Proceedings of the Winter School ``Geometry and Physics''
  (Srn\'\i, 1991)}, number~30, pages 15--20, 1993.

\bibitem{Wal1}
F.~Waldhausen.
\newblock Algebraic {$K$}-theory of topological spaces. {I}.
\newblock In {\em Algebraic and geometric topology (Proc. Sympos. Pure Math.,
  Stanford Univ., Stanford, Calif., 1976), Part 1}, Proc. Sympos. Pure Math.,
  XXXII, pages 35--60. Amer. Math. Soc., Providence, R.I., 1978.

\bibitem{WalM1}
Friedhelm Waldhausen.
\newblock Algebraic {$K$}-theory of spaces, a manifold approach.
\newblock In {\em Current trends in algebraic topology, Part 1 (London, Ont.,
  1981)}, volume~2 of {\em CMS Conf. Proc.}, pages 141--184. Amer. Math. Soc.,
  Providence, R.I., 1982.

\bibitem{WalM2}
Friedhelm Waldhausen.
\newblock Algebraic {$K$}-theory of spaces, concordance, and stable homotopy
  theory.
\newblock In {\em Algebraic topology and algebraic $K$-theory (Princeton, N.J.,
  1983)}, volume 113 of {\em Ann. of Math. Stud.}, pages 392--417. Princeton
  Univ. Press, Princeton, NJ, 1987.

\end{thebibliography}

\end{document}